\documentclass[a4paper]{article}

\usepackage{amssymb}
\usepackage{amsmath}
\usepackage{latexsym}

\def\andy {\mbox{~and~}}
\def\setC{\mathbb{C}}
\def\setN{{\rm I\!N}}
\def\R{{\rm I\!R}}

\def\bbbN{{\rm I\!N}} 
\def\lr{Lachi\`eze-Rey} 
\def\guill {\textquotedblleft ~}

\def\ie {{i.e.}}

\newcommand{\matrixdd}[4]{\left[\begin{array}{cc}#1&#2\\#3&#4\\ \end{array} \right]}
\newcommand{\bin}[2]{\left(  \begin{array}{cc}#1\\#2\\ \end{array}  \right)}

\title{Laplacian eigenmodes for the three sphere }
\author{M. Lachi{\`e}ze-Rey\\
Service d'Astrophysique, C. E. Saclay\\91191 Gif sur Yvette cedex, France}
\begin{document}
\maketitle \abstract{ The vector space ${\cal V} ^{k }$ of the eigenfunctions of the Laplacian on  the three sphere   $S^3$, corresponding to the same eigenvalue  $\lambda _{k} = -k~(k+2)$, has dimension
$(k+1)^2$. After  recalling the standard bases  for ${\cal V} ^{k }$, we introduce a new basis B3, constructed from the reductions to $S^3$ of peculiar homogeneous harmonic polynomia involving null vectors. We give the transformation laws between  this basis  and the usual hyper-spherical harmonics. 

Thanks to the quaternionic representations of $S^3$ and SO(4), we are able to write explicitely the transformation properties of B3, and thus of any   eigenmode, under an arbitrary rotation of SO(4). 
This offers the possibility to select those functions of ${\cal V} ^{k }$  which remain invariant under a chosen rotation of SO(4).  When the rotation is an holonomy transformation of a spherical space $S^3/\Gamma$, this gives a method to calculates the eigenmodes of $S^3/\Gamma$, which remains an open probleme in general. We illustrate our method by (re-)deriving the eigenmodes of lens and prism space. In a forthcoming paper, we present the derivation for dodecahedral space.}

\section{Introduction}

The eigenvalues of the Laplacian $\Delta$ of $S^3$ are of the form $\lambda _{k} = -k~(k+2)$,  where $k \in \bbbN^+$.     For a given value of $k$, they span the  eigenspace 
    ${\cal V} ^{k }$ of   dimension $(k+1)^2$. This vector space constitutes the   $(k+1)^2$ dimensional  irreductible  representation of SO(4), the isometry group of $S^3$.
    
    There are two commonly used bases (hereafter B1 and B2) for   ${\cal V} ^{k }$ which generalize in some sense (see below) the usual spherical harmonics $Y_{\ell  m}$  for the two-sphere.
The functions of these bases have a friendly  behavior under some of the   rotations of SO(4); this generalizes the  property of the  $Y_{\ell  m}$  to be eigenfunctions of the angular momentum operator in $\R ^3$.  However, these functions show no peculiar properties under the {\sl general} rotation  of SO(4). 

Excepted for some cases (lens and prism spaces, see below), the search for  the eigenmodes of the spherical spaces of the form     $S^3/\Gamma$ remains an open problem. Since those are eigenmodes of $S^3$ which remain invariant under the rotations of $\Gamma$, it is clear that this search requires an understanding of the rotation properties of the basis functions under SO(4).

    The task of this paper is to examine  the rotation properties of  the eigenfunctions of     ${\cal V} ^{k }$, as a preparatory work for    the search for eigenfunctions of  $S^3/\Gamma$ (in particular for dodecahedral  space). This will be done through the introduction of  a new basis B3
     of ${\cal V} ^{k }$ (in the case $k$ even), for which the  rotation properties can be explicitely calculated:
    following a new procedure (that was already applied to $S^2$ in \cite{lach})  we generate a system of coherent states on ${\cal V} ^{k }$. We    extract  from it   a   basis B3 of ${\cal V} ^{k }$,  which seems to  have been ignored in the literature, and  presents original properties. Each function  $\Phi ^{k}_{IJ}$
of this basis  B3  is  defined as [the reduction to $S^3$  of] an homogeneous harmonic polynomial in $\R ^4$, which takes the very simple form $(X \cdot N)^{k}$. Here, the dot  product extends   the  Euclidean [scalar] dot  product of $\R ^4$ to its  complexification $\setC ^4$, and $N$ is a null   vector of $\setC ^4$, that we specify below. 
After defining these functions, we show that they form a basis of ${\cal V} ^{k }$, and we give the explicit  transformation formulae between   B2 and B3. 

The    properties of the basis B3 differ from those of  the two other bases, and make it more convenient  for particular applications. In particular, it is possible to calculate explicitely its  rotation properties, under an arbitrary rotation  of  SO(4), by using their quaternionic representation (section \ref{Rotations}). This allows to find those  functions which remain invariant  under an  arbitrary 
rotation. In section  \ref{lensprism}, we apply these result to (re-)derive the eigenmodes of lens and prism space.  

\section{Harmonic functions}
 
 A  function $f$  on $S^3$ is an eigenmode [of the Laplacian]  if it verifies $\Delta f=\lambda f$.
It is known that    eigenvalues are of the form $\lambda _{k}= -k~(k +2)$, $k \in \setN ^{+}$. The corresponding eigenfunctions
generate the eigen[vector]space  ${\cal V}^{k}$, of dimension $(k+1)^2$, which 
  realizes an irreducible unitary representation of the group SO(4). 
  
{\bf First basis}

      I call   B1 the most  widely used   basis for ${\cal V}^{k}$    provided by  the hyperspherical  harmonics \begin{equation} \label{ }
B1 \equiv ({\cal Y}_{k \ell  m}\propto Y_{\ell  m} ),~\ell = 1..k ,~m=-\ell ..\ell .
\end{equation} 
It generalizes the usual spherical  harmonics $Y_{\ell  m}$ on the sphere. In fact,   
it can be shown (\cite{bander}, \cite{erd} p.240,\cite{fry})  that a basis of this type exists on any 
sphere  $S^{n}$. Moreover,      \cite{erd} \cite{fry}  show  that the B1  basis  for $S^n$  is \guill naturally generated " by the  B1 basis  for $S^{n-1}$. In this sense, the B1  basis  for $S^3$ is  generated by the  usual spherical harmonics $Y_{\ell m }$ on the 2-sphere $S^2$.

The generation process involves    harmonic polynomials constructed from  null complex vectors (see below). The basis B1  is in fact  based on the reduction of the representation of SO(4) to   representations of SO(3): each ${\cal Y}_{k \ell  m}$ is an  eigenfunction  of an SO(3) subgroup  of SO(4) which leaves a selected point  of $S^3$ invariant. This make these functions useful when one considers the action of that peculiar SO(3) subgroup. But they show no simple behaviour  under a  general rotation.  We will no more use this basis.

{\bf Second basis}

By group theoretical arguments, \cite{bander} construct a different  ON basis of $V^k$, which is specific to $S^3$:\begin{equation} \label{B2}
B2 \equiv (T_{k;m_1,m_2}) ,~ m_1, m_2=-k/2...k/2, \end{equation} where  $m_1$ and $m_2$ vary  independently    by entire increments (and, thus, take entire or semi-entire values according to the parity of $k$).   In the spirit of the  construction refered above, B2 may be seen as    generated from   a  different choice   of spherical harmonics on $S^2$. The bases B1 and B2 appear respectively adapted to the systems of hyperspherical and toroidal (see below)  coordinates to describe $S^3$.  

The formula (27) of \cite{bander}, reduced to the three-sphere, shows that the elements of this basis take a very convenient form if we use 
 {\sl toroidal coordinates}   (as they are called by  \cite{Leh})
 on the  three sphere $S^3$: $(\chi, \theta, \phi)$ spanning  the  range $0 \le \chi \le \pi/2$, $0 \le \theta \le 2\pi$ and $0 \le \phi \le 2\pi$. 
 They are conveniently defined  (see \cite{Leh} for a more complete description) from an   isometric  embedding of $S^3$  in $\R ^4$ (as the hypersurface $x \in \R ^4; ~\mid x \mid =1$):
 \begin{displaymath}
      \left\{
      \begin{array}{ccc}
          x^0 & = & r~\cos \chi ~\cos \theta \\
          x^1 & = & r~\sin \chi ~\cos \phi   \\
          x^2 & = & r~\sin \chi ~ \sin \phi  \\
 x^3 & = & r~\cos \chi ~\sin \theta \\
     \end{array}
      \right.
\end{displaymath}
 where  $ (x^{\mu}),~\mu=0,1,2,3$,  is a point of $\R ^4$. As shown in   \cite{Leh}),  these coordinates  appear naturally associated to some  isometries.
 
Very simple manipulations show that, with these coordinates, each  eigenfunction of B2 takes the form:
\begin{equation} \label{ }
T_{k;m_1,m_2} (X )\equiv 
t_{k;m_1,m_2}(\chi)~ e^{i\ell \theta}~~ e^{im \phi}, \end{equation} 
where the  $t_{k;m_1,m_2}(\chi)$ are   polynomials in $\cos \chi$ and $\sin \chi$ and we wrote, for simplification, $\ell \equiv  m_1+m_2,~m \equiv m_2-m_1$. 

To have a convenient expression, we report this formula in the harmonic equation expressed in coordinates $\chi,\theta,\phi$. This leads to a second order differential   equation (cf. equ 15 of \cite{Leh}). The solution  is   proportional  to a Jacobi polynomial:
  $t_{k;m_1,m_2} (\chi) \propto  \cos ^{\ell}\chi  ~ \sin ^m \chi   ~P^{m,\ell}_d(\cos  2 \chi),~ d \equiv  k/2-m_2$. Thus, we have the final expression for the basis B2
 \begin{equation}  \label{B2modes}
T_{k;m_1,m_2} (X )= C_{k;m_1,m_2}~  [\cos \chi ~ e^{i \theta}] ^{\ell  } ~[\sin \chi ~ e^{i \phi}] ^{m  }  ~P^{( m,\ell)}_d[\cos (2 \chi)],\end{equation} 
  with   $C_{k;m_1,m_2}~\equiv \frac{\sqrt{ (k+1)}}{\pi} ~\sqrt{\frac{ (k/2+m_2)!~(k/2-m_2)! }{(k/2+m_1)!~Ê(k/2-m_1)!}}$ from normalization requirements (the variation ranges of $m_1$ anf $m_2$ imply that the quantities under factorial sign are entire and positive).
Note also the useful  proportionality relations:
$$      \cos  ^{\ell }\chi   ~ \sin   ^{ m    }\chi  ~P^{( m  ,   \ell)}_{\frac{k-\ell-m)}{2}}(\cos 2 \chi)
     \propto
 \cos  ^{\ell   } \chi ~ \sin  ^{ -m    } \chi  ~P^{( -m  ,   \ell)}_{\frac{k-\ell+m)}{2}}(\cos 2 \chi)
     $$
     $$\propto
 \cos  ^{-\ell   }  \chi~ \sin  ^{ m    } \chi  ~P^{( m  ,  - \ell)}_{\frac{k+\ell-m}{2}}(\cos 2 \chi).$$

   The term $~\zeta ^{m}~\xi ^{n} \equiv e^{i \ell \theta}~e^{i m \phi}$ in (\ref{B2modes}) defines the rotation properties of $T_{k;m_1,m_2} $ under a specific  subgroup of SO(4). This properties generalizes the properties of the spherical harmonics on the two-sphere $S^2$, 
   to be  eigenfunctions  of the rotation operator $P_x$. This advantage has been used 
  by \cite{Leh} to calculate (from a slightly different basis) the eigenmodes of lens or prism spaces (see below section \ref{lensprism}).
   However, the  $T_{k;m_1,m_2} $ have no simple rotation properties under the general rotation of SO(4). This     motivates the search for  a different basis of ${\cal V}^{k}$.
  
Note that the    basis functions $T_{k;m_1,m_2}$  have  also been  introduced in \cite{erd} (p. 253), with their expression in Jacobi Polynomials.  Note also that  they are the complex counterparts of those proposed by \cite{Leh}  (their equ. 19)  to find the eigenmodes of lens and prism spaces. The variation range of the indices  $m_1,m_2$ here (equ. \ref{B2})  is equivalent to 
 their condition \begin{equation} \label{ }
\mid \ell \mid + \mid m \mid  \le k,~\ell +m =k, ~\mod(2),
\end{equation}
through the correspondence $\ell=m_1+m_2,~m=m_2-m_1$.

\subsection{Complex null vectors}

A complex   vector $Z \equiv (Z^0,Z^1,Z^2,Z^3)$ is an  element of $\setC ^4$. We extend the Euclidean scalar product in $\R ^4$  to the complex (non Hermitian) inner product  $Z \cdot Z'\equiv \sum_{\mu} ~Z^{\mu}~Ê(Z^{\prime})^{\mu}$, $\mu=0,1,2,3$. A \emph{null}  vector $N$ is defined as having   zero    norm $N \cdot N \equiv \sum _{\mu} ~N^{\mu} N^{\mu}=0 $  (in which  case, it may be   considered as a point on the isotropic cone in $\setC ^4$).
It is  well known  that     polynomials  of the form $(X \cdot N)^{k}$, homogeneous of degree $k$, are  harmonic if and only if $N$ is a null vector. This results from 
$$\Delta _0 (X \cdot N)^{k} \equiv  \sum _{\mu}  \partial _{\mu} ~\partial _{\mu}~Ê(X \cdot N)^{k}=
k ~\left( Ê\sum _{\mu} (N _{\mu} ~N _{\mu})\right)~Ê(X \cdot N)^{k -1}=0,$$ where $\Delta _0$ is the Laplacian of $\R ^4$.
Thus, the    restrictions of  such polynomials  belong to ${\cal V}^{k}$. As we mentioned above, peculiar null vectors have been used in \cite{erd} and  \cite{fry} to generate the bases B1 and B2. 

To construct a third basis B3, let us first  define a family of null  vectors \begin{equation}
\label{ }N({a,b}) \equiv   (cos a  ,~i~\sin b ,~i~\cos b ,~\sin a ),  
\end{equation}indexed by two   angles  $a$ and $b$ describing the unit circle (they define coherent states in $\R ^4$). 

 The polynomial $[X \cdot N(a,b)]^k$ is harmonic and, thus, can be decomposed onto the  
  basis B2. It is easy to check that, like  the scalar product $X \cdot N(a,b)$, this  polynomial    depends on $a$ and $b$  only through  the combinations $e^{i(\theta- a)}$ and $e^{i(\phi+b)}$, with their conjugates. This    implies that   its decomposition on B2 takes   the form \begin{equation}
\label{devPB2}[X \cdot N(a,b)]^k=\sum _{m_1,m_2} ~P_{k;m_1,m_2}~ T_{k;m_1,m_2}(X)~ e^{-ia(m_1+m_2)}~~ e^{ib(m_2-m_1)}, \end{equation}
where the coefficients $P_{k;m_1,m_2}$ do not depend on $a,b$.
Now we intend to find a basis of ${\cal V}^{k} $ in the form of such polynomials.

\subsection{An new basis}
\subsubsection{Roots of unity}
 
 To do so, we consider the   $(k +1)^{th}$ complex roots of unity  which are the powers $\rho ^I$ of \begin{equation}
\label{}\rho \equiv e^{\frac{2 i \pi }{k+1}} \equiv \cos \alpha +i \sin \alpha  , ~\alpha  \equiv \frac{2 \pi}{k+1}. 
\end{equation}We recall the fundamental property, which will be widely used thereafter:\begin{equation}
\label{property}
\sum_{n=0}^k \rho ^{nI}=(k+1)~\delta ^{Dirac}_I,
\end{equation}where the equallity in the Dirac must be taken $\mod k+1$.

In a given frame, we consider the family of null vectors  \begin{equation}
\label{ }
N_{IJ} \equiv N(I \alpha,J\alpha)=(\cos I  \alpha  ,~i~\sin J\alpha ,~i~\cos J \alpha ,~\sin I   \alpha ), ~ I,J=0..k
\end{equation} 
and  we  define the functions $\Phi^{k}_{IJ}:~\Phi^{k}_{IJ} (X) \equiv (X \cdot N_{IJ} )^{k}$.  
We   report such a  function  in    equ.(\ref{devPB2}) to obtain  its development  in the basis B2. Then we multiply both terms
 by $\rho^{I(m_1+m_2)-J (m_2-m_1)}$. Making the summations over $I,J$ (each varying from $0$ to $k$), and using (\ref{property}), we obtain, in the case where $k$ is even (that we assume hereafter):
 \begin{equation} \label{Tphi}
{\cal  T}_{k;m_1,m_2}= \frac{1}{(k+1)^2 }~\sum _{I,J=0}^k \rho^{I(m_1+m_2)-J (m_2-m_1)}~ \Phi^{k}_{IJ}, \end{equation}
where $~Ê{\cal  T}_{k;m_1,m_2} \equiv ~P_{k;m_1,m_2}~ T_{k;m_1,m_2}$.

This gives the decomposition of  any $T_{k;m_1,m_2}$ (and thus, of any harmonic function)    as a sum of the $(k+1)^2$ polynomials $\Phi^{k}_{IJ}$, providing  the  new  basis of ${\cal V}^{k}$: \begin{equation} \label{ }
B3\equiv (\Phi^{k}_{IJ}), ~ I,J=0..k~ ~ ~ ~ ~  (k~Ê\mbox{even}). 
\end{equation}    The coefficients $P_{k;m_1,m_2}$ involved in the transformation  are calculated in Appendix A. We   obtain easily the reciprocal   formula  expressing the change of basis:\begin{equation}
\label{ }  \Phi^{k}_{IJ} = \sum _{m_1,m_2=-k/2}^{k/2}~ 
{\cal  T}_{k;m_1,m_2}~  \rho^{-I(m_1+m_2)+J (m_2-m_1)}. \end{equation}

\section{Rotations in $\R ^4$}\label{Rotations}

\subsection{Matrix representations}

The isometries of $S^3$ are  the rotations in the embedding space $\R ^4$. In the usual matrix representation, a rotation is represented by a $4*4$ orthogonal matrix $g \in $ SO(4), acting on the 4-vector $(x^{\mu})$ by matrix product.

In the complex matrix representation, a  point (vector) of $\R ^3$  is represented by the
 $2*2$ complex matrix  $$
X \equiv    \matrixdd{W}{iZ}{i\bar{Z}}{\bar{W}}; ~W \equiv x^0+ix^3, ~Z \equiv x^1+ix^2 \in \setC.$$ 
A  rotation $g$  is represented by two complex $2*2$ matrices $(G_L,G_R)$, so that its action takes the form $X \mapsto G_L~X ~G_R$ (matrix product).
The two matrices  $G_L$ and $G_R$ belong to SU(2). Since SU(2) identifies with $S^3$, any matrix $G_L$  or $G_R$ is of the same form than the matrix $X$ above. Since  SU(2) is also the set of unit norm quaternions, there is   a quaternionic representation for the action of SO(4).

\subsection{Quaternionic notations}

Let us note $j_{\mu},~\mu=0,1,2,3$ the basis of quaternions (the $j_{\mu}$ correspond to the usual $1,i,j,k$ but we do not use this notation here). We have $j_0=1$.  A general quaternion is $q=q^{\mu} j_{\mu}=q^0+q^{i} j_{i}$ (with summation convention; the index $i$ takes the values 1,2,3; the index $\mu$ takes the values 0,1,2,3). Its quaternionic  conjugate is $\bar{q} \equiv q^0-q^{i} j_{i}$.
The scalar product is $q_1 \cdot q_2 \equiv (q_1~\bar{q}_2+q_2~\bar{q}_1)/2$, giving  the 
quaternionic norm $\mid q \mid ^2=\frac{q \bar{q}}{2}=\sum _{\mu} (q^{\mu})^2$.

We represent any point $x=(x^{\mu})$  of $\R ^4$ by the quaternion $q_x \equiv x^{\mu}~j_{\mu}$. The points of the (unit) sphere $S^3$ correspond to units quaternions, $\mid q \mid ^2=1$. 
Hereafter, all quaternions will be unitary (if no otherwise indicated). It is easy to see that, using the coordinates above, a point of $S^3$ is represented by the quaternion $\cos \chi ~\dot{\zeta}+ \sin \chi~ \dot{\xi}~ j_1$, 
where we  define dotted quantities, like  $\dot{\zeta} \equiv \cos \theta + j_3~\sin \theta$,  $\dot{\xi} \equiv \cos \phi + j_3~\sin \phi$, as the quaternionic  analogs of the complex numbers $\zeta= \cos \theta + i~\sin \theta$ and $\xi =\cos \phi + i~\sin \phi$,   \ie, with the imaginary $i$ replaced by the quaternion  $j_3$. 

In quaternionic notation, the rotation $g:~  x \mapsto gx$ is represented by a pair of unit quaternions $(Q_L,Q_R)$ such that $q_x \mapsto q_{gx}=Q_L ~q_x ~Q_R$.

{\bf  Complex quaternions, null quaternions}

The null  vectors $N$ introduced above do not belong to $\R ^4$ but to $\setC ^4$. Thus, they cannot be represented by quaternions, but by complex quaternions. Those are defined exactly like  the usual quaternions, but with complex rather than real coefficients. Note that the pure imaginary $i$ does not coincide with any of the $j_{\mu}$, but commutes with all of them. Also, complex conjugation (star)  and quaternionic conjugation (bar)  must be carefully distinguished. Then it is easy to see that the (null) vectors $N_{IJ}$ defined above  correspond to the complex quaternion $n_{IJ} \equiv \dot{\rho}^I +i~ j_2 ~\dot{\rho}^{J}$. Note that  $\mid  n_{IJ}\mid ^2=0$.

In  quaternionic notations, the basis functions are expressed as \begin{equation}
\label{ }
 \Phi_{IJ} (x)=(N_{IJ} \cdot x)^k=<n_{IJ} \cdot  q_X>^k=\left(\frac{n_{IJ}~\bar{q_X}+q_X~\bar{n}_{IJ}}{2}\right)^k.\end{equation}
Quaternionic notations will help us to check how the basis functions are transformed by the rotations of SO(4).

\subsection{Rotations of functions}

To any  rotation      $g$, is associated its action $\mathbf{R}_g$ on functions:
$\mathbf{R}_g: f \mapsto \mathbf{R}_gf;~\mathbf{R}_g f(x)\equiv f(gx)$. 
Let us  apply this   action to the basis functions:\begin{equation}
\label{ }\mathcal{R}_g  \Phi_{IJ} (x)=\Phi_{IJ} (gx)=<n_{IJ} \cdot (Q_L~q_x~Q_R)>^k.
\end{equation} We consider a function on $S^3$ also  as a functions on the set of  unit  quaternions ($q_x$ is the unit quaternion associated    to the point $x$ of $S^3$).
On the other hand, we may  develop this  function on the basis:\begin{equation}
\label{developed}\mathbf{R}_g  \Phi_{IJ} \equiv \sum _{ij=0}^k~ G_{IJ}^{ij}(g)~\Phi_{ij}.\end{equation} The coefficients $G_{IJ}^{ij}(g)$ of the development, that   we intend to calculate, completely encode the action of the rotation $g$ on the basis B3, and thus on  $V^k$.

To proceed , we introduce three auxiliary   complex quaternions: $$ \alpha \equiv 1+i~j_3 ,~  \beta \equiv j_1-i~j_2=  (1-i~j_3)~j_1  \andy  \delta \equiv - j_1-i~j_2.$$ They have zero norm and obey  the properties  $<\alpha \cdot n_{IJ}>=  \rho ^{I},~\\ <\bar{\alpha} \cdot n_{IJ}>=  \rho ^{-I},~<\beta \cdot n_{IJ}>= \rho ^J,~<\delta \cdot n_{IJ}>= \rho ^{-J}$.
Let us now estimate the relation (\ref{developed})  for the specific quaternion    $
\alpha + R~Ê\bar{\alpha} + S~  \beta +T~\delta$, with $R,S,T$  arbitrary real numbers:  
\begin{equation}
\label{ } ({\cal A}+ R~{\cal A}'+S~{\cal B}+T~{\cal D})^k= \sum _{i j} G_{IJ}^{i j}(g)~< (\rho ^i  + R~Ê\rho ^{-i} +S~Ê\rho ^j  + T~Ê\rho ^{-j} )^k,\end{equation} 
where 
${\cal A} \equiv <Q_L~ \alpha ~Q_R \cdot n_{IJ}>$,
${\cal A}' \equiv <Q_L~ \bar{\alpha} ~Q_R \cdot n_{IJ}>$,
${\cal B} \equiv <Q_L~ \beta ~Q_R \cdot n_{IJ}>$,
${\cal D} \equiv <Q_L~ \delta ~Q_R \cdot n_{IJ}>$ characterize the rotation. (Note that these quantities depend on $I$ and $J$).

We develop and identify the powers of the exponents $R,S,T$:
$$ {\cal A}^q~Ê {\cal A}'^{p-q} ~{\cal B}^{r} ~{\cal D}^{k-p-r}= \sum _{ij} G_{IJ}^{ij}(g)~ \rho ^{i(2q-p)} ~Ê\rho ^{j(2r-k+p)}.$$ This holds for  $0\le q \le p,~0\le r \le k- p,~0\le p \le k$. After definition of  the new indices $A \equiv q+r,~B \equiv q-r+k-p$, which  both vary from 0 to $k$, the previous equation takes the form
 $$
  \left(\frac{{\cal A} {\cal B}}{ {\cal A}'   {\cal D}}\right)^{A/2}
   \left(\frac{{\cal A} {\cal D}}{ {\cal A}'   {\cal B}}\right)^{B/2}
 \left(\frac{{\cal A} {\cal A}' }{    {\cal B} {\cal D}}\right)^{p/2}
   \left(\frac{ {\cal A}'    {\cal B} {\cal D}}{{\cal A} }\right)^{k/2}
= \sum _{ij} G_{IJ}^{ij}(g)~ \rho ^{i(A+B-k) +j(A-B)}.$$  
This  holds for any value of $A,B,p$. A consequence is that $ {\cal A} {\cal A}'=    {\cal B} {\cal D}$, which can be checked directly. Finally, 
 $$ {\cal U }^{A}
{\cal V }^{B}~\left(  {\cal A}'    \right)^{k}= \sum _{ij} G_{IJ}^{ij}(g)~ \rho ^{i(A+B-k)} ~Ê\rho ^{j(A-B)},$$  with ${\cal U }\equiv  \left(\frac{  {\cal B}}{ {\cal A}' }\right)$,
${\cal V} \equiv  \left(\frac{{\cal A} }{    {\cal B}}\right)$.

Taking into account the properties of the roots of unity, this equation has the solution
\begin{equation}
\label{G} G_{IJ}^{ij}= \frac{\left(  {\cal A}'    \right)^{k}}{(k+1)^2} \sum _{A,B=0 }^k~ \rho ^{-i(A+B-k)} ~Ê\rho ^{-j(A-B)}~{\cal U }^{A}
~{\cal V }^{B}.\end{equation} When a rotation is specified, there is no difficulty to estimate the associated values of ${\cal A}'  $, ${\cal U }$, ${\cal V}$, and thus of these coefficients which  completely encode the transformation properties of the basis functions of   $V^k$ under SO(4). 

In the next section, we apply these results to rederive the eigenmodes of Lens or Prism space. 
In the next paper \cite{lachdodec}, we   take for $g$ the    generators  of  $\Gamma$, the group  of holonomies of the dodecahedral  space. This will allow   the selection of  the invariant functions, which constitute its  eigenmodes.

\section{Lens and Prism space}\label{lensprism}
The eigen modes for Lens and Prism space have been found by \cite{Leh}. Here we derive them again for  illustration of  our method. 

 \subsection{Lens space}

An holonomy transformation of a lens space takes the form, in complex  notation,  \begin{equation}
\label{ }
 G_L=\matrixdd{e^{i\frac{ \psi_1 +\psi_2}{2}}}{0}{0}{e^{-i\frac{\psi_1 +\psi_2}{2}}},
 ~G_R=\matrixdd{e^{i\frac{\psi_1 -\psi_2}{2}}}{0}{0}{e^{-i\frac{\psi_1 -\psi_2)}{2}}}G.
\end{equation} Its action on a vector of  $ \R ^4$ takes the form  \begin{equation}
\label{ }
X \equiv    \matrixdd{W}{iZ}{i\bar{Z}}{\bar{W}} \mapsto   G_L~ÊX~G_R=Ê \matrixdd{W e^{i \psi _1 }}{iZ e^{i \psi _2}}{i\bar{Z e^{-i \psi _2 }}}{\bar{W}e^{-i \psi _1 }}.
\end{equation} In this simple case,  $W  \equiv x^0+i x^3$, $Z \equiv x^1+i x^2$ are transformed into $W~Êe^{i \psi _1 }$ and $Z~Êe^{i \psi _2}$ respectively. This corresponds to the   quaternionic notation 
  \begin{equation}
\label{ }
Q_L=\dot{w}_1~\dot{w}_2,~Q_R=\dot{w}_1/\dot{w}_2,~\dot{w}_i \equiv \cos (\psi_i/2)+j_3~\sin  (\psi_i/2).
\end{equation}

The rotation  is expressed in the simplest way in the toroidal coordinates, since it acts as $\theta \mapsto \theta +\psi _1,~\phi \mapsto \phi + \psi _2$. From the expression (\ref{B2modes}) of the basis functions (B2), it result their transformation law : 
$$\mathbf{R}_g: T_{k;m_1,m_2} \mapsto T_{k;m_1,m_2} ~Êe^{\ell \psi_1+m \psi _2}.$$ This leads directly  to the invariance condition $\ell \psi_1+m \psi _2=0 \mod 2 \pi$. Using the standard notation for a   lens space $L(p,q)$,  namely
$$\psi _1=2\pi /p,~\psi _2=2\pi ~q /p,$$ we are led to the conclusion:\\

{\framebox{\shortstack{
the  eigenmodes of lens space $L(p,q)$ are all  linear combinations of 
$T_{k;\underline{m_1},\underline{m_2}}$, \\ where the underlining means that the indices verify the condition\\
$\underline{m_1}+\underline{m_2}+q (\underline{m_2}-\underline{m_1})=0,~ \mbox{modulo}(p)$.
}}}

\subsection{Prism space}

The two generators are  single action rotations ($G_R=0$). The first generator, analog  to the lens case above, with $\psi  _1=\psi_2=2\pi /2P$, provides    the first condition $ 
\ell +m=0, ~\mod{2P}$ which takes the form \begin{equation}
\label{ }
\underline{m_2}=0 \mod P .\end{equation} This implies that $k$ must be even.

The second generator has the complex matrix  form  $G=G_L=\matrixdd{0}{-i}{-i}{0}$, which  corresponds to the quaternion $Q_L=Q=-j_1$. Easy calculations lead to 
${\cal A}=- \rho ^J$,
${\cal A}'= \rho ^{-J}$,
${\cal B} =\rho ^I$,
${\cal D} =- \rho ^{-I}$. 
Reporting in (\ref{G}) gives 
\begin{equation}
 G_{IJ}^{ij}= \frac{\rho ^{(i-J)~k}}{(k+1)^2} 
 \sum _{A,B=0 }^k~ 
 \rho ^{A~(-i-j+I+J)+B~Ê(-i+j-I+J)} ~Ê
~(- 1)^{B}.\end{equation}
This formula, together with those expressing the change of basis between B2 and B3, allow to return  to the rotation properties of the basis B2 which take the simple form:\begin{equation}
\label{ }\mathbf{R}: {\cal T}_{k;m_1,m_2} \mapsto (-1)^{m_2+k/2}~{\cal T}_{k;m_1,-m_2}.
\end{equation}It results immediately that the $G$-invariant functions are
combinations of $ {\cal T}_{k;m_1,m_2} + (-1)^{m_2+k/2}~{\cal T}_{k;m_1,-m_2}$. 

Finally, \\

{\framebox{\shortstack{
the eigenfunctions of the Prism space are combinations of \\ $
{\cal T}_{k;m_1,\underline{m_2}} + (-1)^{\underline{m_2}+k/2}~{\cal T}_{k;m_1,-\underline{m_2}}, ~\forall m_1 $;~ $k$ even.   
}}}\\
According to the parity of $k/2$, the functions ${\cal T}_{k;m_1,0}$ are included or not, from which simple counting give the multiplicity as \\
$(k+1)~(1+[k/2P])$, for $k$ even ($[...]$ means entire value),\\
$(k+1)~[k/2P]$, for $k$ odd, in accordance with   \cite{Ikeda}.
  
\section{Conclusion}

We have shown that $ V^k$, the space of eigenfunctions of the Laplacian of $S^3$ with a given eigenvalue $\lambda _k$ ($k$ even) admits   a  new basis B3. In contrary to standard bases (B1 and B2) which  show specific   rotation properties  under selected  subgroups of SO(4),  it is possible  to calculate explicitely  the rotation properties of   B3 under any  rotation of the group SO(4), as well as to calculate  the functions  invariant under  this rotation. This opens the door to the calculation of eigenmodes of spherical space.
 The eigenfunctions of  lens and prism spaces had  been calculated by \cite{Leh}, by using a basis related to B2 (its real, rather than complex, version). We rederived them to illustrate the properties of the bases.

In a subsequent paper \cite{lachdodec}, we apply these results to the search of the eigenfunctions of the dodecahedral space $S^3 /\Gamma $, where $\Gamma=D^*_P $ is the binary dihedral group of order $4P$. Those functions,  still presently  unknown,  are the  eigenfunctions of  $S^3$ which remain invariant under the elements of $\Gamma$.

\subsection{Appendix A}
Let us evaluate  the function \begin{equation}
\label{ }
Z^k_{\ell m} (X)\equiv \sum _{ IJ=0}^k ~\rho ^{\ell I-Jm}  ~{\Phi}_{IJ}^{k} (X)
\end{equation} $$
 =  2^{-k}~\sum _{IJ} ~\rho ^{\ell I -mJ}~\left[ \cos \chi ~(\zeta \rho ^{-I}+\frac{1}{\zeta \rho ^{-I}})  +  \sin \chi  ~(\xi \rho ^J-\frac{1}{\xi \rho ^J})\right]^{k},$$where we defined  $\zeta \equiv e^{i\theta}$ and $\xi  \equiv e^{i\phi}$.  After     development of  the power with the binomial coefficients, the sum becomes  
\begin{eqnarray}\label{fgh}
\sum _{IJ} ~\rho ^{\ell I -mJ}~ \sum _{p=0}^{k} \bin{k}{p}~
[\cos \chi ~(\zeta \rho ^{-I}+\frac{1}{\zeta \rho ^{-I}})]^{k-p}~
[ \sin \chi  ~(\xi \rho ^J-\frac{1}{\xi \rho ^J})]^{p}.\end{eqnarray}

Let us  write the identities
\begin{equation} \label{ }
\rho ^{\ell I}~(\zeta \rho ^{-I}+\frac{1}{\zeta \rho ^{-I}})^{k-p}=\sum _{r=0}^{k-p}~\bin{k-p}{r}~ \zeta ^{2r+p-k}  ~\rho ^{-I(2r+p-k-\ell)},
\end{equation}
\begin{equation}
\label{ }
\rho ^{-mJ}~(\xi \rho ^J-\frac{1}{\xi \rho ^J})^{p}=\sum _{q=0}^{p} \bin{p}{q}   ~\xi ^{2q-p} (-1)^{p-q}   ~\rho ^{J(2q-p-m)},\end{equation} 
that we  insert into (\ref{fgh}). After  summing over $I,J$, and    rearranging the terms, we obtain:
\begin{eqnarray}
 Z_{\ell m}(X) = 2^{-k}~\zeta ^{\ell}~\xi ^{m}~k!  ~ \sum _{q}^{k} ~
\frac{ (-1)^{q-m}~(\cos \chi)^{k-2q+m} ~(\sin \chi)^{2q-m}~}{q!~(q-m)!~(\frac{k+\ell-2q+m}{2})!~ (\frac{k-\ell-2q+m}{2})!~ }.\end{eqnarray}
This formula results from the fact that, through   (\ref{property}), the summations over $I,J$ imply $p=2q-m$ and $2r=\ell+k+m-2q$, that we have reported.  The range of the summation over $q$ is defined by the  conditions\begin{equation}
\label{ }
0 \le \ell+k+m-2q \le 2k+2m-4q \le 2k,~~~0 \le q \le 2q-m \le k.
\end{equation}
 
Rearrangements of the previous formula, inserting $u \equiv \cos(2 \chi) =2\cos^2 \chi-1=1-2\sin ^2 \chi$, lead to \begin{eqnarray}
 Z_{\ell m}(X) =\frac{ 2^{-3k/2}~\zeta ^{\ell}~\xi ^{m}~k! ~(1+u)^{\frac{  \ell  }{2}}~~(1-u)^{\frac{ m  }{2}}}{ (m+d)!~(\ell +d)!}\\
  \nonumber \sum _{q} ~\bin{m+d}{i}~\bin{\ell +d}{d-i}
~(1+u)^{i} ~(u-1)^{d-i} ,\end{eqnarray}
 where we have defined $i \equiv \frac{k+m-  \ell }{2}-q$ and $d \equiv \frac{k-  \ell   -  m  }{2}$. Verification shows that the range defined as above gives exactly the development formula for the Jacobi polynomial. The comparison with (\ref{Tphi})  gives the coefficient  
   $$ P_{k;m_1,m_2}  =\frac{ 2^{-k}~k! ~}{ (k/2- m_1)!~(k/2+  m_1)!~(k+1)^2 ~ C_{k;m_1,m_2}}  $$
   $$ =\frac{ 2^{-k}~~\pi~k! ~(k+1)^{-5/2}}{\sqrt{{ (k/2+m_2)!~(k/2-m_2)!  (k/2+m_1)!~Ê(k/2-m_1)!}}}  $$

\end{document}